%
\magnification=\magstep1   
\input amstex
\UseAMSsymbols
\input pictex
\vsize=23truecm
\NoBlackBoxes
\parindent=20pt
  \font\rmk=cmr8 
\font\ttk=cmtt8
  

\def\mod{\operatorname{mod}}

\def\End{\operatorname{End}}

\def\Ext{\operatorname{Ext}}

   
\def\arr#1#2{\arrow <1.5mm> [0.25,0.75] from #1 to #2}

\centerline{\bf Wild Algebras: Two Examples.}
	\bigskip
\centerline{Claus Michael Ringel}
	\bigskip
Let $k$ be a field and $\Lambda$ a finite-dimensional $k$-algebra (associative, with $1$). 
A recent preprint by Chindris, Kinser and Weyman draws the attention to the present sheer
ignorance concerning the possible behavior of wild algebras. The aim of this note is
to exhibit two examples which answer questions mentioned in the paper.
	\bigskip
{\bf 1\. Wild, Schur-representation-finite algebras of 
global dimension $2$.}
		\medskip 
Chindris-Kinser-Weyman [CKW] have asked whether there do exists wild, Schur-representation-finite 
algebras of finite global dimension. Consider the algebra $\Lambda$ given by the quiver $Q$
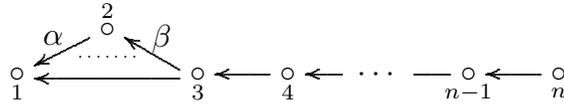
$$
{\beginpicture
\setcoordinatesystem units <1cm,.5cm>
\multiput{$\circ$} at 0 0  1 1   2 0  3 0  5 0  6 0 /
\put{$\cdots$} at 4 0 
\arr{2.8 0}{2.2 0}
\arr{3.6 0}{3.2 0}
\plot 4.8 0  4.4 0 /
\arr{5.8 0}{5.2 0}
\arr{0.8 0.8}{.2 0.2}
\arr{1.8 0.1}{1.2 0.8}
\arr{1.8 -.1}{0.2 -.1}
\put{$\ssize 1$} at 0 -.4
\put{$\ssize 2$} at 1 1.4
\put{$\ssize 3$} at 2 -.4
\put{$\ssize 4$} at 3 -.4
\put{$\ssize n-1$} at 5 -.4
\put{$\ssize n$} at 6 -.4
\put{$\alpha$} at 0.4 0.8
\put{$\beta$} at 1.6 0.8 
\setdots <1mm>
\plot 0.7 0.4  1.3 0.4 /
\endpicture}
$$
with $n\ge 3$ vertices and the zero relation $\alpha\beta$. Since $Q$
is directed, its global dimension of $\Lambda$ is finite. Actually, an easy calculation shows that the 
global dimension is equal to $2$.

The vertex $2$ is a node in the sense of Martinez [M], thus there is a natural bijection between
the indecomposable non-simple representation of $\Lambda$, and the 
indecomposable non-simple representations of the following quiver:
$$
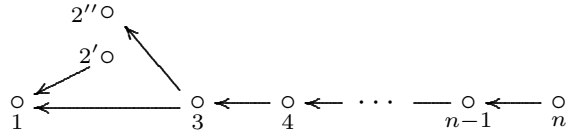
{\beginpicture
\setcoordinatesystem units <1cm,.5cm>
\multiput{$\circ$} at 0 0  1 1  1 2   2 0  3 0  5 0  6 0 /
\put{$\cdots$} at 4 0 
\arr{2.8 0}{2.2 0}
\arr{3.6 0}{3.2 0}
\plot 4.8 0  4.4 0 /
\arr{5.8 0}{5.2 0}
\arr{0.8 0.8}{.2 0.2}
\arr{1.8 0.3}{1.2 1.7}
\arr{1.8 -.1}{0.2 -.1}
\put{$\ssize 1$} at 0 -.4
\put{$\ssize 2'$} at 0.8 1.1
\put{$\ssize 2''$} at 0.75 2
\put{$\ssize 3$} at 2 -.4
\put{$\ssize 4$} at 3 -.4
\put{$\ssize n-1$} at 5 -.4
\put{$\ssize n$} at 6 -.4
\endpicture}
$$
For $n\ge 9$, this is a wild quiver, thus, for $n\ge 9$, the algebra $\Lambda$ is wild.
	\medskip
{\bf Lemma.} {\it If $M$ is an indecomposable
representation of $\Lambda$ such that $M_\alpha, M_\beta$ both are
non-zero, then there is a non-zero endomorphism $\phi$ of $M$ with $\phi^2 = 0.$}
	\medskip
Proof. If $M_\alpha \neq 0,$ then the simple module $S(2)$ is a factor module of $M$,
if $M_\beta\neq 0$, then $S(2)$ is a submodule of $M$. Thus, if both $M_\alpha, M_\beta$ are
non-zero,  we obtain an endomorphism $\phi$
of $M$ with image $S(2)$. We must have $\phi^2 = 0$, since otherwise $S(2)$ is a direct
summand of $M$; but then $M = S(2)$, impossible.
	\medskip
{\bf Corollary.} {\it If $M$ is a representation whose endomorphism ring $\End M$ is a division ring, then
$M_\alpha= 0$ or $M_\beta = 0,$ thus $M$ is a representation of the $\Bbb D_n$ quiver obtained from $Q$ by deleting $\alpha$,
or a representation of the $\Bbb A_n$ quiver obtained by deleting $\beta$.}
	\medskip
This shows that there are only finitely many isomorphism classes of representations $M$ such that 
$\End(M)$ is a division ring, thus $\Lambda$ is Schur-representation-finite. 
	\bigskip\bigskip 
{\bf 2\. A strictly wild algebra without a wild tilted factor algebra.}
		\medskip
A conjecture by Yang Han [H] quoted in [CKW] asserts that any 
strictly wild algebra should have a factor algebra which is a wild tilted algebra. 

Let $k$ be an algebraically closed field and $\Lambda$ the $k$-algebra
with quiver and relations 
$$
{\beginpicture
\setcoordinatesystem units <3cm,1cm>
\multiput{$\circ$} at 0.1 0  .9 0 /
\arr{0.8 .2}{0.2 .2}
\arr{0.8 -.0}{0.2 -.0}
\setquadratic
\plot 0.2 -0.3  0.5 -.5  0.8 -.3 /
\arr{0.78 -.32}{0.8 -.3}
\put{$0$} at 0.1 -.3
\put{$1$} at 0.9 -.3
\put{$\alpha_0$} at 0.5 0.4
\put{$\alpha_1$} at 0.5 -.2
\put{$\beta$} at 0.5 -.8
\put{$\beta\alpha_i\beta = 0$ \ for \ $i=0,1$.} [l] at 1.5 0 
\endpicture}
$$ 
Note that the subquiver given by the arrows $\alpha_0,\alpha_1$ is the Kronecker quiver.
The representations of the Kronecker quiver are called Kronecker modules.
Since we assume that $k$ is algebraically closed, a Kronecker module is
simple regular if and only if it is indecomposable and of length 2, and the
simple regular Kronecker modules $R(\lambda)$
are indexed by the elements $\lambda\in \Bbb P^1(k) = k\cup\{\infty\}$. 

Let $\Cal C$ be the full subcategory of $\mod \Lambda$ given by all $\Lambda$-modules
$M$ such the restriction of $M$ to the Kronecker quiver is a direct sum of simple regular
Kronecker modules. Clearly, $\Cal C$ is
an abelian category (with an exact embedding into $\mod\Lambda$). 
If we endow the simple regular Kronecker module 
$R(\lambda)$ with $\beta$ as zero map, we obtain a simple object in $\Cal C$,
we denote it again by $R(\lambda)$.
Now let $\Cal D$ be the full subcategory of $\Cal C$ consisting of all modules 
$M$ in $\Cal C$
which have a submodule $M'$ which is a direct sum of copies of $R(\infty)$ such that
$M/M'$ is a direct sum of modules of the form $R(\lambda)$ with $\lambda\in k.$
       \medskip
{\it The category $\Cal D$ is an abelian category with exact embedding functor into
$\mod\Lambda$; its simple objects are the modules $R(\lambda)$ with $\lambda\in 
\Bbb P^1(k)$, and we have $\Ext^1_{\Cal D}(R(\lambda),R(\mu)) = k$ provided
$\lambda \in k$ and $\mu = \infty$ and equal to zero, otherwise.}
For $\lambda\in k$, a non-trivial element of 
$\Ext^1_{\Cal D}(R(\lambda),R(\omega)) = k$ is given by using $\beta.$ 
	\medskip 
Thus, the quiver $\Delta$ of $\Cal D$ is a subspace quiver, the number of sources in $\Delta$ is equal to 
$|k|$.
This shows that $\Cal D$ and therefore $\mod\Lambda$ is strictly wild.
Of course, instead of taking such a large subcategory $\Cal D$, it would be sufficient to
ask in the definition of $\Cal D$ that $M/M'$ 
is a direct sum of modules of the form $R(\lambda)$ with $\lambda$ belonging to {\bf a fixed 
5-element subset} of $k$ (so that one obtains the 5-subspace quiver).
	  \medskip
On the other hand, any factor algebra of $\Lambda$ is given by a quiver with at most 2
vertices (and some relations), and a wild tilted algebra with at most 2 vertices
is hereditary (it is a generalized Kronecker algebra with at least 3 arrows). Of course, $\Lambda$ has no such factor
algebra. 
	\bigskip
{\bf References}
	\medskip
\item{[CKW]} Chindris, C\.,  Kinser, R\., Weyman, J.: Module varieties and representation type of finite-dimensional
   algebras.  arXiv:1201.6422
\item{[H]} Han, Y.: Controlled wild algebras.
   Proc. London Math. Soc. (3), 83 (2001), 279-298.
\item{[M]} Martinez-Villa, R.: Algebras stably equivalent to 1-hereditary. In:
  Representation Theory II. Springer Lecture Notes in Math. 832 (1980), 396-431.  
\bigskip

{\rmk
\noindent 
Shanghai Jiao Tong University, Shanghai 200240, P. R. China, and \newline 
  King Abdulaziz University, PO Box 80200,  Jeddah, Saudi Arabia. \newline
E-mail: 
{\ttk ringel\@math.uni-bielefeld.de}
}

\bye